\documentstyle{article}
\begin{document}
\title{Witten-Helffer-Sj\"{o}strand Theory for a Generalized Morse
Function}
\author{Hon-kit Wai}
\maketitle
\begin{abstract} In this paper, we extend the Witten-Helffer-Sj\"{o}strand
theory from Morse functions to
generalized Morse functions. In this case, the spectrum of the Witten
deformed Laplacian $\Delta(t)$, for
large t, can be seperated into the small eigenvalues (which tend to 0
as $t\rightarrow\infty$), large and
very large eigenvalues (both of which tend to $\infty$ as
$t\rightarrow\infty$). The subcomplex
$\Omega_{0}^{\ast}(M,t)$ spanned by eigenforms corresponding
to the small and large eigenvalues of $\Delta(t)$
is finite dimensional. Under some mild conditions, it is shown
that $(\Omega_{0}^{\ast}(M,t),d(t))$ converges
to a geometric complex associated to the generalized Morse
function as $t\rightarrow\infty$.
\end{abstract}
\section{Introduction and Statement of Results}

The purpose of this paper is to extend the
Witten-Helffer-Sj\"{o}strand theory
(cf.[W],[H-S]) for a Morse
function on compact manifold to a generalized Morse function. Such a
generalized Morse function has all
critical points either non-degenerate or of birth-death type, i.e. in
some neighbourhood of the critical
point and with respect to a convenient coordinate system, the function
can be written

\[f(x_{1},x_{2},\ldots,x_{n}) = f(0)-x_{1}^{2}-\ldots-x_{k}^{2}+
x_{k+1}^{2}+\ldots+x_{n-1}^{2}+ax_{n}^{3} \]
for some $a \not= 0$.

The interest of generalized Morse function comes from the following
theorem due to H.Chaltin(cf.[I]).

\begin{quote}
{\bf Theorem}: If $\pi: E \rightarrow B$ is a smooth  bundle with
fibre a compact manifold M, then there exists ]
$f: E \rightarrow R$ so that for any $t \in B$

\[f_{t} = f\mid_{\pi^{-1}(t)}: M_{t} = \pi^{-1}(t) \rightarrow R\] is
a generalized Morse function.
\end{quote}
It is easy to see that in general one cannot have such a statement
with $f_{t}$ a Morse function.

Now, let us state the results of this paper.  Let the eigenvalues of
the operator \[-\frac{d^{2}}{dx^{2}}+9x^{4}-6x \]
be \[0<e_{1}<e_{2}\leq \ldots \leq e_{l}\leq \ldots  \]
(See Lemma 2 in \S 2 for proof.)

 Suppose $M^{n}$ is a compact orientable Riemannian manifold, f be a
 generalized Morse function on M.

Suppose $x_{1}^{k},\ldots,x_{m_{k}}^{k}$ are all the non-degenerate
critical points of f, of index k, $y_{1}^{k},
\ldots,y_{m_{k}^{\prime}}^{k}$, are all the critical points of
birth-death type, of index k.  Also, let $a_{j}^{(k)}
\epsilon R$ be associated with $y_{j}^{k}$ so that in some neighbourhood
of $y_{j}^{k}$ and with respect to a
suitable oriented co-ordinate system,

\[f(x_{1},\ldots,x_{n}) = f(0)-x_{1}^{2}-\ldots-x_{k}^{2}+x_{k+1}^{2}+
\ldots+x_{n-1}^{2}+a_{j}^{(k)}x_{n}^{3}\]

Suppose for simplicity that
\begin{equation}
 \mid a_{1}^{(0)} \mid < \mid a_{2}^{(0)} \mid <\ldots< \mid
 a_{m_{0}^{\prime}}^{(0)} \mid < \mid a_{1}^{(1)}
 \mid <\ldots< \mid a_{m_{1}^{\prime}}^{(1)} \mid <\ldots< \mid
 a_{m_{n-1}^{\prime}}^{(n-1)} \mid
\end{equation}
(in fact, the Witten-Helffer-Sj\"{o}strand theory is very similar with
minor modifations without assuming (1))

  Also, let g be a Riemannian metric on M so that in the above
  co-ordinate system near the critical points $ x_{j}^{k}$
  or $y_{j}^{k}$, g is the canonical metric on $R^{n}$.

Consider the Witten deformation of the de Rham complex
$(\Omega^{\ast}(M),d(t))$ with \[d(t) = e^{-tf}de^{tf}:
\Omega^{\ast}(M) \rightarrow \Omega^{\ast}(M)\]
Consider the deformed Laplacian
 \[ \Delta(t) = d(t)d^{\ast}(t)+d^{\ast}(t)d(t)\]

When the above canonical coordinates near the critical points are used,

\begin{equation}
 \Delta(t) = \Delta+t^{2}{\mid df\mid}^{2}+tA
\end{equation}
where \[A = \sum_{i,j=1}^{n}\frac{\partial^{2}f}{\partial x_{i}\partial
x_{j}}[dx_{i},i_{\partial_{j}}]\]
and $i_{\partial_{j}}$ denotes the contraction along the vector field
$\partial_{j}$, $dx_{i}$ is the exterior
multiplication by the form $dx_{i}$ and $[dx_{i},i_{\partial_{j}}]$
denotes the commutator
$dx_{i}i_{\partial_{j}}-i_{\partial_{j}}dx_{i}$.

There are two cases.

{\bf Case 1}: $x_{j}^{k}$ is non-degenerate.

\begin{equation}
f(x) = f(0)-x_{1}^{2}-\ldots-x_{k}^{2}+x_{k+1}^{2}+\ldots+x_{n}^{2}
\end{equation}
 \[{\mid df \mid}^{2} =  4(x_{1}^{2}+\ldots+x_{n}^{2})\]
\begin{equation}
A = -2\sum_{i=1}^{k}[dx_{i},i_{\partial_{i}}]+2\sum_{i=k+1}^{n}[dx_{i},
i_{\partial_{i}}]
\end{equation}
{\bf Case 2}: $y_{j}^{k}$ is of birth-death type.

\begin{equation}
f(x) = f(0)-x_{1}^{2}-\ldots-x_{k}^{2}+x_{k+1}^{2}+\ldots+x_{n-1}^{2}+
a_{j}^{(k)}x_{n}^{3}
\end{equation}
\[{\mid df \mid}^{2} = 4(x_{1}^{2}+\ldots+x_{n-1}^{2})+9(a_{j}^{(k)})^{2}
x_{n}^{4}\]

\begin{equation}
A = -2\sum_{i=1}^{k}[dx_{i},i_{\partial_{i}}]+2\sum_{i=k+1}^{n-1}[dx_{i},
i_{\partial_{i}}]+6a_{j}^{(k)}x_{n}[dx_{n},
i_{\partial_{n}}]
\end{equation}
For each critical point c, define the $'localized'$ operator
$\overline{\Delta}_{c}(t): C^{\infty}(\Lambda^{\ast}(R^{n}))
\rightarrow C^{\infty}(\Lambda^{\ast}(R^{n}))$ which is given by (2)
where $A$ is (4) if $c=x_{j}^{k}$ is nondegenerate,
respectively (6) if $c=y_{j}^{k}$ is birth-death. The Laplace operator
in (2) then is the canonical Laplace operator
corresponding the canonical metric in $R^{n}$.  The operator
$\overline{\Delta}_{c}(t)$ then extends uniquely to a
self-adjoint positive unbounded operator in $L^{2}(\Lambda^{\ast}(R^{n}))$.

Now suppose $\overline{\Delta}(t)$ is the $'localized'$ operator associated
to a critical point of birth-death type.
Since $L^{2}(\Lambda^{\ast}(R^{n})) \cong L^{2}(\Lambda^{\ast}(R^{n-1}))
\otimes L^{2}(\Lambda^{\ast}(R))$,
$\overline{\Delta}(t)$ can be written as

\begin{eqnarray*}
 \overline{\Delta}(t) & = &   \{\Delta_{R^{n-1}} +4t^{2}(x_{1}^{2}+\ldots+
 x_{n-1}^{2})+2t\sum_{i=1}^{n-1}
 \epsilon_{i}[dx_{i},i_{\partial_{i}}]\} \otimes id \\
    &    & \mbox{}  + id \otimes \{\Delta_{R}+9(a_{j}^{(k)})^{2}t^{2}
	x_{n}^{4}+6a_{j}^{(k)}tx_{n}[dx_{n},
	i_{\partial_{n}}]\}  \end{eqnarray*}
where  \[\epsilon_{i} = \left\{ \begin{array}{rl}
                                  -1   &  \mbox{if $ 1\leq i\leq k$}     \\
                                   1   &  \mbox{if $ k+1\leq i\leq n-1$}
                                \end{array}
                        \right. \]
and $\Delta_{R^{n-1}}$, $\Delta_{R}$ are the Laplace operators on $R^{n-1}$
and $R$ respectively.

By Corollary in \S 2, $\Delta_{R}+9(a_{j}^{(k)})^{2}t^{2}x_{n}^{4}+
6a_{j}^{(k)}x_{n}[dx_{n},i_{\partial_{n}}]:
L^{2}(\Lambda^{\ast}(R)) \rightarrow L^{2}(\Lambda^{\ast}(R))$ has
discrete spectrum with eigenvalues

\[0<e_{1}(\mid a_{j}^{(k)}\mid t)^{2/3}<e_{2}(\mid a_{j}^{(k)}\mid t)^{2/3}
\leq e_{3}(\mid a_{j}^{(k)}\mid t)^{2/3}
\leq \ldots \]

Each eigenvalue has a multiplicity of 2 with corresponding eigenvectors
consisting of a 0-form and a 1-form.

Let $\Delta^{k}(t) = \Delta(t)\mid_{L^{2}(\Lambda^{k}(M))}$.

Let $0 \leq E_{1}(t) \leq E_{2}(t) \leq \ldots \leq E_{l}(t) \leq \ldots $
be all the eigenvalues of $\Delta^{k}(t)$.

Suppose for simplicity that
\[e_{1}\mid a_{m_{n-1}^{'}}^{(n-1)}\mid^{2/3}<e_{2}\mid a_{1}^{(0)}
\mid^{2/3}\]
This together with (1) imply that $e_{1}\mid a_{m_{k}^{'}}^{(k)}
\mid^{2/3}<e_{2}\mid a_{1}^{(k-1)}\mid^{2/3}$ for all k.
\newtheorem{theorem}{Theorem}
\begin{theorem}[Quasi-classical limit of eigenvalues]

\[\lim_{t \rightarrow \infty}E_{1}(t) = \ldots = \lim_{t \rightarrow
\infty}E_{m_{k}}(t) = 0\]
\[\lim_{t \rightarrow \infty}\frac{E_{m_{k}+1}(t)}{t^{2/3}} = e_{1}(\mid
a_{1}^{(k-1)}\mid)^{2/3} < \lim_{t \rightarrow
\infty}\frac{E_{m_{k}+2}(t)}{t^{2/3}} = e_{1}(\mid a_{2}^{(k-1)}
\mid)^{2/3} < \ldots \]
\[\ldots <\lim_{t \rightarrow \infty}\frac{E_{m_{k}+m_{k-1}^{'}+
m_{k}^{'}}(t)}{t^{2/3}} = e_{1}(\mid
a_{m_{k}^{'}}^{(k)}\mid)^{2/3} \]
\[ \left( < \lim_{t \rightarrow \infty}\frac{E_{m_{k}+m_{k-1}^{'}+
m_{k}^{'}+1}(t)}{t^{2/3}} = e_{2}(\mid a_{1}^{(k-1)}
\mid)^{2/3} < \ldots  \right) \]
\end{theorem}
{\bf Remarks}:1.  In fact, the eigenvectors corresponding to $E_{1}(t),
\ldots, E_{m_{k}}(t)$  are localized at the
non-degenerate critical points of index k, while the eigenvectors
corresponding to $E_{m_{k}+1}(t),\ldots,E_{m_{k}+
m_{k-1}^{'}+m_{k}^{'}}(t)$ are localized at the birth-death critical
points of index k-1 and k. However, the
eigenvectors are not necessarily localized at a single critical point.

2.  If all the critical points of f are non-degenerate, then the above
theorem should be formulated as follows
(cf. [S] p219):

{\bf Theorem}: \[\lim_{t \rightarrow \infty}E_{1}(t) = \ldots = \lim_{t
\rightarrow \infty}E_{m_{k}}(t) = 0\]
\[0 < \lim_{t \rightarrow \infty}\frac{E_{m_{k}+1}(t)}{t} \leq \lim_{t
\rightarrow \infty}\frac{E_{m_{k}+2}(t)}{t}
\leq \ldots \]

\noindent Let us index $a_{1}^{(0)},\ldots,a_{m_{0}^{'}}^{(0)},a_{1}^{(1)},
\ldots,a_{m_{n-1}^{'}}^{(n-1)}$ by $b_{1},
\ldots,b_{N}$ where $N=\sum_{k=0}^{n-1}m_{k}^{'}$.

Also, for $0\leq l\leq n-1, 1\leq j\leq m_{l}^{'}$ let
\[I_{j}^{(l)}(\epsilon)=[e_{1}(\mid a_{j}^{(l)}\mid)^{2/3}-\epsilon,e_{1}
(\mid a_{j}^{(l)}\mid)^{2/3}+\epsilon]\]
Choose an $\epsilon$ small enough so that the family of intervals
\[\left\{[0,\epsilon],I_{j}^{(l)}(\epsilon),[e_{2}(\mid a_{1}^{(0)}
\mid)^{2/3}-\epsilon,\infty)\right\}_{0\leq l\leq
n-1,1\leq j\leq m_{l}^{'}}\]
is pairwise disjoint.  The pairwise disjointness is satisfied if
$\epsilon$ is a positive number smaller than
\[min_{i}\left(\frac{e_{1}}{2}\mid a_{1}^{(0)}\mid^{2/3},
\frac{e_{1}}{2}(\mid b_{i+1}\mid^{2/3}-\mid
b_{i}\mid^{2/3}),\frac{1}{2}(e_{2}\mid a_{1}^{(0)}\mid^{2/3}-e_{1}
\mid a_{m_{n-1}^{'}}^{(n-1)}\mid^{2/3})\right) \]

As a consequence of Theorem 1, for t sufficiently large and $\epsilon$
satisfying the above disjointness condition,
we have
\[Spec(t^{-2/3}\Delta^{k}(t)) \subset [0,\epsilon]\cup \left(
\cup_{l=k-1,k;1\leq j\leq m_{l}^{'}}I_{j}^{(l)}(\epsilon)
\right)\cup[e_{2}(\mid a_{1}^{(0)}\mid )^{2/3}-\epsilon,\infty) \]
and \[\left\{ \begin{array}{l}
               Card\left(Spec(t^{-2/3}\Delta^{k}(t))\cap[0,\epsilon]
			   \right) = m_{k}  \\
               Card\left(Spec(t^{-2/3}\Delta^{k}(t))\cap I_{j}^{(l)}
			   (\epsilon) \right) = 1 \\
               (for l=k-1,k;1\leq j\leq m_{l}^{'})
              \end{array}
      \right. \]

Define \[\Omega_{small}^{k}(M,t) = Span\{\psi(t)\in L^{2}(\Lambda^{k}(M))
|\Delta^{k}(t)\psi(t)=E(t)\psi(t),
t^{-2/3}E(t)\in[0,\epsilon]\} \]
\begin{eqnarray*}
 \Omega_{large,l,j}^{k}(M,t)& = &Span\{\psi(t)\in L^{2}(\Lambda^{k}(M))|
 \Delta^{k}(t)\psi(t)=E(t)\psi(t), \\
                            &   &t^{-2/3}E(t)\in [e_{1}(\mid a_{j}^{(l)}
							\mid)^{2/3}-\epsilon,e_{1}(\mid a_{j}^{(l)}
							\mid)^{2/3}+\epsilon]\} \\
 \Omega_{v.large}^{k}(M,t)  & = &Span\{\psi(t)\in L^{2}(\Lambda^{k}(M))|
 \Delta^{k}(t)\psi(t)=E(t)\psi(t), \\
                            &   &t^{-2/3}E(t)\in [e_{2}-\epsilon,\infty)\}
\end{eqnarray*}
\[\left(\Omega_{0}^{\ast}(M,t),d(t)\right)=
\left(\Omega_{small}^{\ast}(M,t),d(t)\right)\perp\left(\perp_{k,j}
\left(\Omega_{large,k,j}^{\ast}(M,t),d(t)\right)\right) \]
\newtheorem{corollary}{Corollary}
\begin{corollary}
$\left(\Omega^{\ast}(M),d(t)\right)$ is equal to
\[\left(\Omega_{small}^{\ast}(M,t),
d(t)\right)\perp\left(\perp_{k,j}
\left(\Omega_{large,k,j}^{\ast}(M,t),d(t)\right)\right)\perp
\left(\Omega_{v.large}^{\ast}(M,t),d(t)\right) \]
$\left(\Omega_{small}^{\ast}(M,t),d(t)\right),\left(\Omega_{0}^{\ast}
(M,t),d(t)\right)$ are complexes of finite
dimensional vector spaces which calculate the de Rham cohomology of M.
\end{corollary}
{\bf Remark}:  Observe that $\left(\Omega_{large,k,j}^{\ast}(M,t),d(t)
\right)$ has dimension 2.  It is spanned by a
k-form and a (k+1)-form localized at $y_{j}^{k}$ (the localization of
the forms at $y_{j}^{k}$ is due to (1)).

As in the Helffer-Sj\"{o}strand theory for a generic pair (f,g),
$\left(\Omega_{0}^{\ast}(M,t),d(t)\right)$ converges as
$t\rightarrow\infty$ to a geometric complex, which can be described as
follows.

Let f be a self-indexing generalized Morse function, i.e.

\[\left\{ \begin{array}{ll}
                   f(x_{j}^{k})=k      & \mbox{if $x_{j}^{k}$ is a
				   non-degenerate critical point of index k}  \\
                   f(y_{j}^{k}) \in (k,k+1)  & \mbox{if $y_{j}^{k}$ is
				   a birth-death critical point of index k}
                  \end{array}
          \right. \]
Let $W_{x_{j}^{k}}^{-}$ be the descending manifold of a non-degenerate
critical point $x_{j}^{k}$.  For a birth-death
critical point $y_{j}^{k}$, choose an open neighbourhood $U_{y_{j}^{k}}$
and a suitable co-ordinate s.t.
\[f(x)=f(0)-x_{1}^{2}- \ldots -x_{k}^{2}+x_{k+1}^{2}+ \ldots +x_{n-1}^{2}+
a_{j}^{(k)}x_{n}^{3} \]
Let \begin{eqnarray*}
     W_{y_{j}^{k}}^{-,0} & = & \mbox{$ \{x\in M|\gamma_{x}(t)\in
	 U_{y_{j}^{k}}\cap R^{k}$ for some t$ \in R $}  \\
                       &   & \mbox{where $\gamma_{x}$ is the trajectory
					   of Grad f s.t.$ \gamma_{x}(0)=x$ \} } \\
     W_{y_{j}^{k}}^{-,1} & = & \mbox{$ \{x\in M|\gamma_{x}(t)\in
	 U_{y_{j}^{k}}\cap (R^{k}\times R_{-})$ for some t$ \in R $} \\
                       &   & where R^{k}\times R_{-}=\{(x_{1}, \ldots ,
					   x_{k},0, \ldots ,0,x_{n})\in R^{n}|x_{n}<0 \} \}
    \end{eqnarray*}

when $a>0$ with obvious modifications when $a<0$.  Then
$W_{y_{j}^{k}}^{-,0},W_{y_{j}^{k}}^{-,1}$ are manifolds
diffeomorphic to $R^{k},R^{k+1}$ respectively. Note that
$W_{y_{j}^{k}}^{-,0}\cap W_{y_{j}^{k}}^{-,1}=\emptyset$

Define the descending manifold
\[W_{y_{j}^{k}}^{-}=W_{y_{j}^{k}}^{-,0}\cup W_{y_{j}^{k}}^{-,1} \]
which is then a manifold with boundary diffeomorphic to $R_{+}^{k+1}$.
The ascending manifold $W_{y_{j}^{k}}^{+}$ is
defined similarly.

Suppose the ascending and descending manifolds for any two critical
points intersect transversally, then
$\left\{W_{x_{j}^{k}}^{-},W_{y_{j}^{k}}^{-,0},W_{y_{j}^{k}}^{-,1}
\right\}$ form a CW-complex (see $\S3$ for more details).
While the incidence number between $W_{x_{j}^{k}}^{-}$ and
$W_{x_{i}^{k+1}}^{-}$ is given by the intersection number
between the ascending and decending manifolds in $f^{-1}(k+
\frac{1}{2})$ i.e. between
$W_{x_{j}^{k}}^{+}\cap f^{-1}(k+\frac{1}{2})$ and
$W_{x_{i}^{k+1}}^{-}\cap f^{-1}(k+\frac{1}{2})$, the incidence number
is 1 between $W_{y_{j}^{k}}^{-,0}$ and $W_{y_{j}^{k}}^{-,1}$.  However,
those between $W_{x_{j}^{k}}^{-}$ and
$W_{y_{j}^{k}}^{-,i}$ may be non-trivial (i=0,1).  Let us denote the
above described chain complex by $(C_{\ast}(M,f),
\partial)$ (with $C_{k}(M,f)=Span\{W_{x_{j}^{k}}^{-},W_{y_{j}^{k}}^{-,0},
W_{y_{j}^{k-1}}^{-,1}\})$, its dual cochain
complex by $(C^{\ast}(M,f),\delta)$.

Also, let us rescale the complex $(\Omega_{0}^{\ast}(M,t),d(t))$ to be
\[(\Omega_{0}^{\ast}(M,t),\tilde{d}(t)) = \left(\Omega_{small}^{\ast}(M,t),
e^{t}\sqrt{\frac{\pi}{2t}}d(t) \right)\perp
\left(\perp_{k,j}\left(\Omega_{large,j,k}^{\ast}(M,t),d(t) \right)
\right) \]

\begin{theorem}
There exists $f^{\ast}(t):(\Omega_{0}^{\ast}(M,t),\tilde{d}(t))
\rightarrow
(C^{\ast}(M,f),\delta)$ which is a morphism
of co-chain complexes s.t.
\[f^{\ast}(t) = I +O(t^{-1})\]
w.r.t. some suitably chosen bases.
\end{theorem}

\noindent{\bf Definitions}:  1.(i). Suppose $x_{j}^{k},x_{i}^{k+1}$ are
two non-degenerate critical points, $\gamma$
be a generalized trajectory betwen $x_{i}^{k+1}$ and $x_{j}^{k}$, i.e.
$\gamma$ is a piecewise smooth curve with
singularities at the birth-death points $y_{1},\ldots,y_{n(\gamma)}$ and
$$\gamma=\gamma_{x_{i}^{k+1}y_{1}}\cup \{y_{1}\}\cup \gamma_{y_{1}y_{2}}
\cup \{y_{2}\}\cup \ldots \cup
\gamma_{y_{n(\gamma)}x_{j}^{k}}$$
where $\gamma_{y_{l}y_{l+1}}$ is a trajectory between $y_{l}$ and
$y_{l+1}$.
Then one can associate $\epsilon=\pm 1$
to the trajectories $\gamma_{x_{i}^{k+1}y_{1}},\gamma_{y_{l}y_{l+1}},
\gamma_{y_{n(\gamma)}x_{j}^{k}}$ as in the
Witten-Morse theory.  Then define
$$\epsilon_{\gamma}^{new}=(-1)^{n(\gamma)}
\epsilon_{\gamma_{x_{i}^{k+1}y_{1}}}\left(\prod_{l=1}^{n(\gamma)-1}
\epsilon_{\gamma_{y_{l}y_{l+1}}}\right)
\epsilon_{\gamma_{y_{n(\gamma)}x_{j}^{k}}}$$
(ii) Suppose $x_{j}^{k}$ is an non-degenerate critical point and
$y_{i}^{k}$ a birth-death critical point, $\gamma$
be a generalized trajectory between them.  With the above notation for
$\gamma$ and $y_{1}=y_{i}^{k}$, define
$$\epsilon_{\gamma}^{new}=(-1)^{n(\gamma)}\left(\prod_{l=1}^{n(\gamma)-1}
\epsilon_{\gamma_{y_{l}y_{l+1}}} \right)
\epsilon_{y_{n(\gamma)}x_{j}^{k}}$$
2. (i) The (generalized) incidence number between two critical points
is defined as follows:
$$I(x_{i}^{k+1},x_{j}^{k})=\sum_{\gamma}\epsilon_{\gamma}^{new}$$
$$I(y_{i}^{k},x_{j}^{k})=\sum_{\gamma}\epsilon_{\gamma}^{new}$$
where $\gamma$ is a generalized trajectory between the initial and
end point.

(ii) Here we recall that the (ordinary) incidence number between two
critical points(non-degenerate or birth-death) is

\[i(x_{i}^{k+1},x_{j}^{k})=\sum_{\gamma}\epsilon_{\gamma}\]
\[i(y_{i}^{k},x_{j}^{k})=\sum_{\gamma}\epsilon_{\gamma}\]
where $\gamma$ is a trajectory between the two critical points.

{\bf Remark}: Observe that in general
$$\epsilon_{\gamma}^{new}\not =\epsilon_{\gamma}$$
for a trajectory between two critical points.  For example, if $\gamma$
is a trajectory between a birth-death point
$y_{i}^{k}$ and a non-degenerate critical point $x_{j}^{k}$, then
$\epsilon_{\gamma}^{new}=-\epsilon_{\gamma}$

With the above definition, we can reformulate Theorem 2 as follows:

\noindent{\bf Theorem $2'$: (Helffer-Sj\"{o}strand)} There exist
orthonormal bases $\left\{E_{x_{j}^{k}}(t)\right\}$ of
$\Omega_{small}^{k}(M,t)$, $\left\{E_{y_{j}^{k}}^{0}(t),
E_{y_{j}^{k}}^{1}(t)\right\}$ of
$\Omega_{large,k,j}^{\ast}(M,t)$ s.t.
\[<E_{x_{i}^{k+1}}(t),d(t)E_{x_{j}^{k}}(t)> = e^{-t}\left(\sqrt
{\frac{t}{\pi}}\sum_{\gamma}\epsilon_{\gamma}^{new} +
O(t^{-1/2}) \right) \]
\[<E_{y_{j}^{k}}^{1}(t),d(t)E_{y_{j}^{k}}^{0}(t)> = \sqrt{e_{1}}
(a_{j}^{k})^{1/3}t^{1/3} + O(t^{1/6}) \]
\[<E_{y_{j_{1}}^{k}}^{i_{1}}(t),d(t)E_{y_{j_{2}}^{l}}^{i_{2}}(t)>=0
\mbox{ if $j_{1}\not=j_{2}$ for t sufficiently
large } \]
\[<E_{x_{i}^{k}}(t),d(t)E_{y_{j}^{l}}^{i^{'}}(t)>=
<E_{y_{j}^{l}}^{i^{'}}(t),d(t)E_{x_{i}^{k}}(t)>=0 \mbox{ for t
sufficiently large } \]

where $\sum_{\gamma}\epsilon_{\gamma}^{new}=I(x_{i}^{k+1},x_{j}^{k})$
is the incidence number between $x_{j}^{k}$
and $x_{i}^{k+1}$ defined above.
\[  \]
Inside the complex $(C^{\ast}(M,f),\delta)$, there is a subcomplex
$(C_{nd}^{\ast}(M,f),\delta)$ such that
$$dimC_{nd}^{k}(M,f)=m_{k}$$
where $m_{k}$ is the number of non-degenerate critical points of index k.
Note that this subcomplex is not generated
by the non-degenerate critical points, since the latter in general do
not generate a subcomplex.  Instead the
subcomplex is obtained by applying Lemma 3 repeatedly as is done in \S 3.
See \S 3 for details.

\noindent{\bf Theorem $2''$}: $f^{\ast}(t)\mid_{\Omega_{small}^{\ast}
(M,t)}: \left(\Omega_{small}^{\ast}(M,t),
\tilde{d}(t)\right) \longrightarrow (C^{\ast}(M,f),\delta)$ is an
injective homomorphism of co-chain complexes whose
image complex converges to $(C_{nd}^{\ast}(M,f),\delta)$ in
$(C^{\ast}(M,f),\delta)$ as $t\rightarrow\infty$, more
precisely,
\[f^{k}(t)\left(E_{x_{j}^{k}}(t)\right)=\hat{e}_{x_{j}^{k}}+O(t^{-1})
\mbox{ in } C^{\ast}(M,f) \]
where $\hat{e}_{x_{j}^{k}}=e_{x_{j}^{k}}+\sum_{l}I(y_{l}^{k},x_{j}^{k})
e_{y_{l}^{k}}^{0}\in C_{nd}^{k}(M,f)$.

Also, using similar consideration,one can extend the result to any
representation $\rho: \pi_{1}(M) \rightarrow GL(V)$
(cf.[BZ]) or any representation $\rho: \pi_{1}(M)\rightarrow
GL(W_{A})$ where $W_{A}$ is finite type Hilbert module
over a finite von Neuman algebra $A$ (cf.[BFKM]).

The above question concerning the extension of Witten-Helffer-
Sj\"{o}strand theory for generalized Morse functions
was raised in Dan Burghelea's course on $L_{2}$-topology.  I would
like to thank him for the problem and help in
accomplishing this work. A parametrized version of the above theory
will be presented in future work in collaboration
with D. Burghelea.

\section{Witten Deformation for a Generalized Morse Function}

Let f be a generalized Morse function on $M^{n}$, y be a critical
point of birth-death type.  Let $(U_{y},\varphi)$
be a chart s.t. $y\in U_{y}$, and
\[f(\varphi^{-1}x)=f(0)-x_{1}^{2}- \ldots -x_{k}^{2}+x_{k+1}^{2}+
\ldots +x_{n-1}^{2}+ax_{n}^{3} \]
Let g be a Riemannian metric on M s.t. $(\varphi^{-1})^{\ast}
(g)=\delta_{ij}$.

Define \[d(t)=e^{-tf}de^{tf} \]
\[\Delta(t)=d(t)d^{\ast}(t)+d^{\ast}(t)d(t) \]
Then in the coordinate system $(U_{y},\varphi)$,
\[\Delta(t)=\Delta+t^{2}\mid df \mid^{2}+2t\sum_{i=1}^{n-1}
\epsilon_{i}[dx_{i},i_{\partial_{i}}]+6ax_{n}t[dx_{n},
i_{\partial_{n}}] \]
where \[\mid df \mid^{2}=4(x_{1}^{2}+ \ldots +x_{n-1}^{2})+9a^{2}
x_{n}^{4} \]

\[\epsilon_{i}=\left\{ \begin{array}{rl}
                        -1 &   \mbox{if $1\leq i\leq k$}     \\
                         1 &   \mbox{if $k+1\leq i\leq n-1$}
                       \end{array}
               \right.  \]
Define $\overline{\Delta}(t):L^{2}(\Lambda^{\ast}(R^{n})) \rightarrow
L^{2}(\Lambda^{\ast}(R^{n}))$ to be given
exactly by the above expression.  Recall that
\begin{eqnarray*}
 \overline{\Delta}(t) & = & \{\Delta_{R^{n-1}}+4t^{2}(x_{1}^{2}+ \ldots
 +x_{n-1}^{2})+2t\sum_{i=1}^{n-1}
 \epsilon_{i}[dx_{i},i_{\partial_{i}}] \} \otimes id    \\
                      &   & id \otimes \{\Delta_{R}+9t^{2}a^{2}x_{n}^{4}+
					  6atx_{n}[dx_{n},i_{\partial_{n}}] \}
\end{eqnarray*}
Observe that the first term is exactly the Witten deformed Laplacian on
$L^{2}(\Lambda^{\ast}(R^{n-1}))$ for the
classical Morse theory (cf.[S],[W]) and that the two operators in
parenthesis commute with each other. Hence to
study the spectrum of $\overline{\Delta}(t)$, it suffices to find out
the spectrum of $\Delta+9t^{2}a^{2}x^{4}+
6tax[dx,i_{\frac{d}{dx}}]$ on $L^{2}(\Lambda^{\ast}(R))$.  Note that
         \[    [dx,i_{\frac{d}{dx}}](f)=-f \]
         \[    [dx,i_{\frac{d}{dx}}](fdx)=fdx  \]
Therefore
\[\left\{\Delta+9t^{2}a^{2}x^{4}+6atx[dx,i_{\frac{d}{dx}}] \right\}(f)=
\left( -\frac{d^{2}}{dx^{2}}+9t^{2}a^{2}x^{4}-
6atx \right)(f) \]
\[\left\{\Delta+9t^{2}a^{2}x^{4}+6atx[dx,i_{\frac{d}{dx}}] \right\}(fdx)=
\left( -\frac{d^{2}}{dx^{2}}+9t^{2}a^{2}x^{4}+
6atx \right)(f)dx  \]
Define $R:L^{2}(R) \rightarrow L^{2}(R)$
\[ \left(Rf \right)(x)=f(-x)\]
Then \[ R^{-1}\left(-\frac{d^{2}}{dx^{2}}+9t^{2}a^{2}x^{4}+6atx
\right)R=-\frac{d^{2}}{dx^{2}}+9t^{2}a^{2}x^{4}-6atx \]
Hence, it is sufficient to consider \[P(at)= -\frac{d^{2}}{dx^{2}}+
9t^{2}a^{2}x^{4}-6atx \]
where $P(t)\equiv -\frac{d^{2}}{dx^{2}}+9t^{2}x^{4}-6tx $ for $t\in R$..

Define $U(\lambda):L^{2}(R) \rightarrow L^{2}(R)$, $\lambda>0$
\[\left(U(\lambda)f \right)(x)= \lambda^{1/2}f(\lambda x) \]
\newtheorem{lemma}{Lemma}
\begin{lemma} For $t>0$,
\[P(t)=U(t^{1/3})\left(t^{2/3}P(1) \right)U(t^{-1/3}) \]
\end{lemma}
\begin{lemma}
(i) $P(1)$ has compact resolvent, hence has discrete spectrum.
\[0\leq e_{1}\leq e_{2} \leq \ldots \leq e_{l}\leq \ldots  \]
(ii) The smallest eigenvalue of $P(1)$ is strictly positive and is simple,
i.e.
\[0<e_{1}<e_{2}\leq \ldots  \]
(iii)Let $\Xi_{1}$ be a normalized eigenfunction of $P(1)$ corresponding
to the smallest eigenvalue $e_{1}$. Then one
can choose $\Xi_{1}$ so that $\Xi_{1}(x)>0$ for all $x\in R$. In
particular $\Xi_{1}(0)^{-1}$ exists.

\end{lemma}

\noindent{\bf Corollary}:  $P(at)$ has spectrum
\[0< e_{1}(\mid at\mid)^{2/3}<e_{2}(\mid at\mid)^{2/3}\leq e_{3}(\mid
at\mid)^{2/3}\leq \ldots \leq e_{l}(\mid at
\mid)^{2/3}\leq \ldots \]

\noindent{\bf Proof of Lemma 2}: (i) P(1) has compact resolvent because
$V(x)=9x^{4}-6x \rightarrow \infty$ as $\mid x
\mid \rightarrow \infty$ and is bounded from below. (cf. [RS] p249)  It
is positive because it is the restriction of
the deformed Laplacian associated with the function $x^{3}$ on the
invariant subspace $L^{2}(R)$.

(ii) Let $\left(-\frac{d^{2}}{dx^{2}}+9x^{4}-6x \right)(f)=0$

Then \[f=c_{1}e^{-x^{3}}+c_{2}e^{-x^{3}}\int_{0}^{x}e^{2u^{3}}du \]
One checks that $f\in L^{2}(R)$ iff $c_{1}=c_{2}=0$. This shows $e_{1}>0$.

To show that $e_{1}$ is simple, one apply the Feynman-Kac formula to
show that $e^{-P(1)}$ has a strictly positive kernel.
This implies that the largest eigenvalue of $e^{-P(1)}$ is simple,
hence the simplicity of the smallest eigenvalue of
$P(1)$.(cf. [GJ] \S3.3)

(iii) The application of Feynman-Kac formula in (ii) shows that
$\Xi_{1}$ can be chosen s.t. $\Xi_{1}(x)>0$ almost
everywhere.  One shows that actually $\Xi_{1}(x)>0$ for all $x\in R$.
$\Box$

Now, we have shown that the smallest eigenvalue of $P(at)$ is
$e_{1}(\mid at\mid)^{2/3}>0$ and is simple.  Since
$-\frac{d^{2}}{dx^{2}}+9a^{2}t^{2}x^{4}+6atx$ is conjugated to
$P(at)$ by an isometry R, the same is true for its
smallest eigenvalue.  Hence $-\frac{d^{2}}{dx^{2}}+9a^{2}t^{2}x^{4}+
6atx[dx,i_{\frac{d}{dx}}]$ on
$L^{2}(\Lambda^{\ast}(R))$ has smallest eigenvalue $e_{1}(\mid at
\mid)^{2/3}$ of multiplicity 2, the corresponding
eigenvectors are $\Xi_{1}(x) \left(\in \Omega^{0}(R) \right)$ and
$\Xi_{1}(-x)dx \left(\in \Omega^{1}(R) \right)$.

Returning to the $'localized'$ operator,
\begin{eqnarray*}
 \overline{\Delta}(t) & = & \left\{\Delta_{R^{n-1}}+4t^{2}(x_{1}^{2}+
 \ldots +x_{n-1}^{2})+2t\sum_{i=1}^{n-1}
 \epsilon_{i}[dx_{i},i_{\partial_{i}}] \right\}\otimes id   \\
                      &   & + id\otimes \left\{\Delta_{R}+9a^{2}
					  t^{2}x_{n}^{4}+6atx_{n}[dx_{n},i_{\partial_{n}}]
					  \right\}
\end{eqnarray*}

The smallest eigenvalue is also $e_{1}(\mid at\mid)^{2/3}$ and of
multiplicity 2, whose eigenvectors are spanned by a
k-form and a (k+1)-form.
\[\omega_{k}(t)=t^{(n-1)/4+1/6}\xi_{1}(t^{1/2}x_{1}) \ldots \xi_{1}
(t^{1/2}x_{n-1})\Xi_{1}(t^{1/3}x_{n})dx_{1}\wedge
\ldots \wedge dx_{k} \]
\[\omega_{k+1}(t)=t^{(n-1)/4+1/6}\xi_{1}(t^{1/2}x_{1}) \ldots \xi_{1}
(t^{1/2}x_{n-1})\Xi_{1}(-t^{1/3}x_{n})dx_{1}\wedge
\ldots \wedge dx_{k}\wedge dx_{n}  \]
where $\xi_{1}(x)$ is the $ground state$ of $-\frac{d^{2}}{dx^{2}}+
4x^{2}$ i.e. $\xi_{1}(x)=e^{-tx^{2}}$.

With the above observations, the proof of Theorem 1 follows essentially
the arguments in [S](pp 219-222). See Appendix
for sketch of proof.

\section{Helffer-Sj\"{o}strand Theory for a Generalized Morse Function}
\newtheorem{proposition}{Proposition}

{\bf Definition}: A pair (f,g) is said to satisfy the Morse-Smale
condition where f is a generalized Morse function if
for any two critical points x and y, the ascending manifold $W_{x}^{+}$
and the descending manifold $W_{y}^{-}$, w.r.t.
$- Grad_{g}f$, intersect transversally.

In the case of a birth-death critical point $y_{j}^{k}$ of index k such
that (5) holds, define
\begin{eqnarray*}
 W_{y_{j}^{k}}^{+,0} & = & \mbox{$\{x\in M|\gamma_{x}(t)\in U_{y_{j}^{k}}
 \cap R^{n-k-1}$ for some $t\in R$}   \\
                   &   & \mbox{where $R^{n-k-1}=\{(0,\ldots,0,x_{k+1},
				   \ldots,x_{n-1},0)\in R^{n}\} \}$}  \\
 W_{y_{j}^{k}}^{+,1} & = & \mbox{$\{x\in M|\gamma_{x}(t)\in U_{y_{j}^{k}}
 \cap(R^{n-k-1}\times R_{+}$ for some $t\in R$}  \\
                   &   & \mbox{where $R^{n-k-1}\times R_{+}=\{(0,
				   \ldots,0,x_{k+1},\ldots,x_{n})\in R^{n}|x_{n}> 0\}
				   \}$}
\end{eqnarray*}
while $W_{y_{j}^{k}}^{-,0},W_{y_{j}^{k}}^{-,1}$ are defined similarly as
in \S 1. Then the ascending and descending
manilfolds are defined as follows:
\[W_{y_{j}^{k}}^{+}=W_{y_{j}^{k}}^{+,0}\cup W_{y_{j}^{k}}^{+,1} \]
\[W_{y_{j}^{k}}^{-}=W_{y_{j}^{k}}^{-,0}\cup W_{y_{j}^{k}}^{-,1} \]

\begin{proposition}  For any pair (f,g), there is a $C^{1}$
approximation $g'$ such that $g=g'$ in a neighbourhood
of the critical points of f and (f,g$'$) satisfies the Morse-Smale
condition.
\end{proposition}
{\bf Proof}:  The proof is the same as in [Sm].

\noindent{\bf Definition}: Let f be a generalized Morse function.
f is said to be self-indexing if
\[ \left\{  \begin{array}{ll}
             f(x_{j}^{k})=k     &  \mbox{if $x_{j}^{k}$ is a
			 non-degenerate critical point of index k}  \\
             f(y_{j}^{k})\in (k,k+1)  &  \mbox{if $y_{j}^{k}$ is
			 birth-death critical point of index k}
            \end{array}
   \right.  \]

\begin{proposition}  For any generalized Morse function f, there
exists a self-indexing generalized Morse function
$f'$ such that f and $f'$ have the same critical points and
corresponding indexes.
\end{proposition}
{\bf Proof}:  The proof is similar as in [M] \S4.

\noindent{\bf Definition}: A pair (f,g) is called a generalized
triangulation if

(i) f is a self-indexing generalized Morse function on M and in a
neighbourhood $U_{c}$ of any critical point c,
one can introduce local coordinates s.t. $g=\delta_{ij}$ and
   \begin{quote}
      (a) if c is non-degenerate
          \[f(x)=f(0)-x_{1}^{2}-\ldots -x_{k}^{2}+x_{k+1}^{2}+
		  \ldots +x_{n}^{2} \]
      (b) if c is of birth-death type
          \[f(x)=f(0)-x_{1}^{2}-\ldots -x_{k}^{2}+x_{k+1}^{2}+
		  \ldots +x_{n-1}^{2}+ax_{n}^{3} \]
   \end{quote}
(ii) (f,g) satisty the Morse-Smale condition.

Let f be a generalized Morse function, $W_{x_{j}^{k}}^{-}$ be the
descending manifold of a non-degenerate critical point.
For a birth-death critical point $y_{j}^{k}$, recall
\begin{eqnarray*}
     W_{y_{j}^{k}}^{-,0} & = & \mbox{$\{x\in M|\gamma_{x}(t)\in
	 U_{y_{j}^{k}}\cap R^{k}$ for some $t\in R$}  \\
                       &   & \mbox{where $\gamma_{x}$ is the
					   trajectory of Grad f s.t. $\gamma_{x}(0)=x$
					   \} }  \\
     W_{y_{j}^{k}}^{-,1} & = & \mbox{$\{x\in M|\gamma_{x}(t)\in
	 U_{y_{j}^{k}}\cap (R^{k}\times R_{-})$ for some $t\in R$}  \\
                       &   & where R^{k}\times R_{-}=\{(x_{1},
					   \ldots ,x_{k},0, \ldots ,0,x_{n})\in R^{n}|
					   x_{n}<0 \} \}
    \end{eqnarray*}

where $R^{k}=\{(x_{1},\ldots,x_{k},0,\ldots,0)\in R^{n}\}$.

Then we have

{\bf Theorem}: Suppose (f,g) is a generalized triangulation, then

(i) $\left\{W_{x_{j}^{k}}^{-},W_{y_{j}^{k}}^{-,0},W_{y_{j}^{k}}^{-,1}
\right\}_{0\leq k\leq n;1\leq j\leq m_{k}
\mbox{ or }m_{k}^{'}}$ is a CW-complex.

(ii)  Let $(C_{\ast}(M,f),\partial)$ be the cellular chain complex of
the above CW-complex (as described in \S1),
$(C^{\ast}(M,f),\delta)$ be its dual co-chain complex.

Then $Int: (\Omega^{\ast}(M),d) \rightarrow (C^{\ast}(M,f),\delta)$
\[\omega \longmapsto \int_{W}\omega \]
is a morphism of co-chain complexes.

\noindent{\bf Proof}:  A proof of this theorem in the case of a Morse
function can be found in [L].  The same argument
also works in the case of a generalized Morse function.  However, a
better argument for this is the following:

(i) One first verify that the partition $\left\{W_{x_{j}^{k}}^{-},
W_{y_{j}^{k}}^{-,0},W_{y_{j}^{k}}^{-,1} \right\} $
is a stratification in the sense of Whitney (see [V2] for definition).
Using the tubular neighbourhood theorem
(Proposition 2.6 [V1]) and the fact that each stratum is diffeomorphic

with an Euclidean space, one concludes that
this partition is a CW-complex.

(ii) The fact that integration is well defined and represents a morphism
of cochain complexes follows from Stokes theorem
in the framework of integration theory  on stratified sets (cf.[F],[V2]).
$\Box$

As a consequence, the composition
\[(\Omega_{0}^{\ast}(M,f,t),d(t)) \stackrel{e^{tf}}{\longrightarrow}
(\Omega^{\ast}(M),d) \stackrel{Int}{\longrightarrow}
(C^{\ast}(M,f),\delta)  \]
is also a morphism of co-chain complexes.

Let $$M_{x_{j}^{k}}=M\setminus \left[\left(\cup_{l\not =j}B(x_{l}^{k},
\eta)\right)\cup \left(\cup_{l}B(y_{l}^{k},
\eta)\right)\right] $$
Let $\Delta_{M_{x_{j}^{k}}}(t)$ be the corresponding Laplace operator on
$M_{x_{j}^{k}}$ with Dirichlet boundary
condition, $\Psi_{x_{j}^{k}}(t)$ be an eigenvector corresponding to the
smallest eigenvalue of
$\Delta_{M_{x_{j}^{k}}}^{k}(t)$ of norm one.

Similarly, let $$M_{y_{j}^{k}}=M\setminus \left[\left(\cup_{l}B(x_{l}^{k},
\eta)\right)\cup \left(\cup_{l\not =j}
B(y_{l}^{k},\eta)\right)\right]$$
With $\Delta_{M_{y_{j}^{k}}}(t)$ similarly defined, let
$\Psi_{y_{j}^{k}}^{0}(t)$, respectively
$\Psi_{y_{j}^{k}}^{1}(t)$ be the smallest eigenvector of
$\Delta_{M_{y_{j}^{k}}}^{k}(t),
\Delta_{M_{y_{j}^{k}}}^{k+1}(t)$ of norm one.

Define $J_{k}(t):C^{k}(M,f) \longrightarrow \Omega^{k}(M)$ by
\[J_{k}(t)\left(e_{x_{j}^{k}} \right)=\Psi_{x_{j}^{k}}(t)  \]
\[J_{k+i}(t)\left(e_{y_{j}^{k}}^{i} \right)=\Psi_{y_{j}^{k}}^{i}(t),
\mbox{ $i=0,1$} \]
where $\left\{e_{x_{j}^{k}},e_{y_{j}^{k}}^{i} \right\}$ is the dual
basis of $\left\{W_{x_{j}^{k}},
W_{y_{j}^{k}}^{i} \right\}$.

Let $Q_{k,small}(t)$, $Q_{k+i,k,j}(t)$ be the orthogonal projection
onto $\Omega_{small}^{k}(M,t)$ and
$\Omega_{large,k,j}^{k+i}(M,t)$ respectively.

Define $Q_{k}(t):J_{k}(t)\left(C^{k}(M,f)\right)\longrightarrow
\Omega_{0}^{k}(M,t)$ by
\[Q_{k}(t)\left(\Psi_{x_{j}^{k}}(t) \right)=Q_{k,small}(t)\left
(\Psi_{x_{j}^{k}}(t) \right) \]
\[Q_{k+i}(t)\left(\Psi_{y_{j}^{k}}^{i}(t) \right)=Q_{k+i,k,j}(t)
\left(\Psi_{y_{j}^{k}}^{i}(t)\right) \]
Let \[H_{k}(t)=\left(Q_{k}(t)J_{k}(t)\right)^{\ast}\left(Q_{k}(t)
J_{k}(t)\right) \]
\[\tilde{J_{k}}(t)=Q_{k}(t)J_{k}(t)(H_{k}(t))^{-1/2} \]
Then $\tilde{J_{k}}(t):C^{k}(M,f) \rightarrow \Omega_{0}^{k}(M,t)$
is an isometry.

Define $E_{x_{j}^{k}}(t)=\tilde{J_{k}}(t)\left(e_{x_{j}^{k}}\right)$,
$E_{y_{j}^{k}}^{i}(t)=\tilde{J_{k}}(t)
\left(e_{y_{j}^{k}}^{i}\right)$.

Note that $E_{x_{j}^{k}}(t)\in \Omega_{small}^{k}(M,t)$,
$E_{y_{j}^{k}}^{i}(t)\in \Omega_{large,k,j}^{k+i}(M,t)$.

\begin{proposition}
\[E_{x_{j}^{k}}(t)=(\frac{2t}{\pi})^{n/4}e^{-tx^{2}}\left(dx_{1}\wedge
\ldots \wedge dx_{k}+O(t^{-1}) \right)  \]
\[E_{y_{j}^{k}}^{0}(t)=(\frac{2t}{\pi})^{\frac{n-1}{4}}e^{-t(x_{1}^{2}+
\ldots
+x_{n-1}^{2})}(\mid at\mid)^{1/6}\Xi_{1}((at)^{1/3}x_{n})\left(dx_{1}
\wedge \ldots \wedge dx_{k}+O(t^{-1})\right) \]
on $U_{x_{j}^{k}}$ and $U_{y_{j}^{k}}$ respectively.
\end{proposition}

\noindent{\bf Remark}:  Note that $(\frac{2t}{\pi})^{\frac{n-1}{4}}$
and $\mid at\mid^{1/6}$ are the normalization
constants for $e^{-t(x_{1}^{2}+\ldots+x_{n-1}^{2})}$ and
$\Xi_{1}((at)^{1/3}x_{n})$ respectively, i.e.
\[\|(\frac{2t}{\pi})^{\frac{n-1}{4}}e^{-t(x_{1}^{2}+\ldots+
x_{n-1}^{2})}\|=\|\mid at\mid^{1/6}\Xi_{1}((at)^{1/3}x)\|=1\]

\noindent{\bf Proof}:  One can follow the argument in [HS] or [BZ].
Note that the term
$(\mid at\mid)^{-1/6}\Xi_{1}((at)^{1/3}x)$ is the ground state of
$-\frac{d^{2}}{dx^{2}}+9a^{2}t^{2}x^{4}-6atx$.
So it is also the first term of the asymptotic expansion.  $\Box$

Recall that we have defined $\epsilon_{\gamma}^{new}$ for a generalized
trajectory between two critical points and
the incidence number $I(x,y)$ between two critical points.  See \S 1
for definitions.  Also define
$Int_{k}=Int\mid_{\Omega^{k}(M)}$. With these definitions, we have

\begin{proposition}
(i)\[Int_{k}e^{tf}\left(E_{x_{j}^{k}}(t)\right)=
(\frac{2t}{\pi})^{\frac{n-2k}{4}}e^{tk}\left(e_{x_{j}^{k}}+
\sum_{l}I(y_{l}^{k},x_{j}^{k})e_{y_{l}^{k}}^{0}+O(t^{-1})\right) \]
(ii)\[Int_{k}e^{tf}\left(E_{y_{j}^{k}}^{0}(t)\right)=
(\frac{2t}{\pi})^{\frac{n-1-2k}{4}}e^{f(y_{j}^{k})}\Xi_{1}(0)
\mid a_{j}^{(k)}t\mid^{1/6}\left(e_{y_{j}^{k}}^{0}+ \sum_{l\not = j}
\beta_{lj}(t)e_{y_{l}^{k}}^{0}+O(t^{-1})\right) \]
(iii)\[Int_{k+1}e^{tf}\left(E_{y_{j}^{k}}^{1}(t)\right)=
(\frac{2t}{\pi})^{\frac{n-1-2k}{4}}e^{tf(y_{j}^{k})}
\frac{\Xi_{1}(0)\mid a_{j}^{(k)}t\mid^{1/6}}{\sqrt{e_{1}}\mid
a_{j}^{(k)}t\mid^{1/3}}\left(\delta\left(e_{y_{j}^{k}}^{0}+
\sum_{l\not = j}\beta_{lj}(t)e_{y_{l}^{k}}^{0}\right)+O(t^{-1})
\right)  \]
\end{proposition}

\noindent{\bf Proof}: We introduce the following notations. Let
$y_{1}^{k},\ldots,y_{m_{k^{'}}}^{k}$ be  all the
birth-death critical points of index k.
Let \[\begin{array}{c}
f(y_{1}^{k})=\ldots=f(y_{r_{1}}^{k})<f(y_{r_{1}+1}^{k})=\ldots=
f(y_{r_{1}+r_{2}}^{k})<f(y_{r_{1}+r_{2}+1}^{k})<\ldots<\ldots \\
<f(y_{r_{1}+\ldots r_{l_{k}-1}+1}^{k})=\ldots=f(y_{r_{1}+\ldots+
r_{l_{k}}}^{k})\end{array} \]
where $m_{k^{'}}=r_{1}+\ldots+r_{l_{k}}$.

Also, for $1\leq q\leq l_{k}$ let
\[L_{q}^{(k)}=\{y_{l}^{k}\mid r_{1}+\ldots+r_{q-1}+1\leq l\leq r_{1}+
\ldots+r_{q-1}+r_{q}\}\]

(i) It is clear from [HS],[BZ] that
\[\int_{W_{x_{j}^{k}}^{-}}e^{tf}E_{x_{j}^{k}}(t)=
(\frac{2t}{\pi})^{\frac{n-2k}{4}}e^{tk}\left(1+O(t^{-1})\right)\]
Suppose $y_{l}^{k}\in L_{1}^{(k)}$, let
\[\partial W_{y_{l}^{k}}^{-,1}=W_{y_{l}^{k}}^{-,0}+i(y_{l}^{k},
x_{j}^{k})W_{x_{j}^{k}}^{-}+\left(\sum_{i\not =j}
i(y_{l}^{k},x_{i}^{k})W_{x_{i}^{k}}^{-}+\sum_{i}i(y_{l}^{k},
y_{i}^{k-1})W_{y_{i}^{k-1}}^{-,1}\right)\]
where $i(y_{l}^{k},x_{j}^{k})$ is the (ordinary) incidence number
between $y_{l}^{k}$ and $x_{j}^{k}$ defined in \S 1.
Denote the expression inside the parenthesis by R, the remainder term.

Therefore, \[\int_{\partial W_{y_{l}^{k}}^{-,1}}e^{tf}E_{x_{j}^{k}}(t)=
\int_{W_{y_{l}^{k}}^{-,0}}e^{tf}E_{x_{j}^{k}}(t)+
i(y_{l}^{k},x_{j}^{k})\int_{W_{x_{j}^{k}}^{-}}e^{tf}E_{x_{j}^{k}}(t)+
\int_{R}e^{tf}E_{x_{j}^{k}}(t)\]
But by Stoke's Theorem,
\[\int_{\partial W_{y_{l}^{k}}^{-,1}}e^{tf}E_{x_{j}^{k}}(t)=
\int_{W_{y_{l}^{k}}^{-,1}}e^{tf}\left(d(t)E_{x_{j}^{k}}(t)
\right)=\int_{W_{y_{l}^{k}}^{-,1}}e^{tf}\left(\sum_{i}\lambda_{i}(t)
E_{x_{i}^{k+1}}(t)\right)\]
for some exponentially decaying functions $\lambda_{i}(t)$.
Since $\mid E_{x_{i}^{k+1}}(t)(x)\mid$ decreases as
$e^{-t\mid f(x)-f(x_{i}^{k+1})\mid}$, $\int_{\partial
W_{y_{l}^{k}}^{-,1}}e^{tf}E_{x_{j}^{k}}(t)$ is of smaller order
compared with $\int_{W_{x_{j}^{k}}^{-}}e^{tf}E_{x_{j}^{k}}(t)$.
The same is true for $\int_{R}e^{tf}E_{x_{j}^{k}}(t)$.

Hence, \[\begin{array}{lll}
\int_{W_{y_{l}^{k}}^{-,0}}e^{tf}E_{x_{j}^{k}}(t) & = &
-i(y_{l}^{k},x_{j}^{k})\int_{W_{x_{j}^{k}}^{-}}e^{tf}
E_{x_{j}^{k}}(t)+O(t^{-1}) \\
        & = & I(y_{l}^{k},x_{j}^{k})(\frac{2t}{\pi})^{\frac{n-2k}{4}}
		e^{tk}+O(t^{-1})
\end{array} \]
One can show by using finite induction on $q$ that for any $y_{l}^{k}
\in L_{q}^{(k)}$
\[\int_{W_{y_{l}^{k}}^{-,0}}e^{tf}E_{x_{j}^{k}}(t)=
I(y_{l}^{k},x_{j}^{k})(\frac{2t}{\pi})^{\frac{n-2k}{4}}e^{tk}+O(t^{-1})\]
This proves (i).

(ii) A direct computation shows that
\[\int_{W_{y_{j}^{k}}^{-,0}}e^{tf}E_{y_{j}^{k}}^{0}(t)=
(\frac{2t}{\pi})^{\frac{n-1-2k}{4}}e^{tf(y_{j}^{k})}\Xi_{1}(0)
\mid a_{j}^{(k)}t\mid^{1/6}\left(1+O(t^{-1})\right) \]
Next note that by choosing a coordinate system $(x_{1}^{(j)},\ldots,
x_{k}^{(j)})$ on $W_{y_{j}^{k}}^{-,0}$ and
extending it to a neighbourhood of $W_{y_{j}^{k}}^{+}=
W_{y_{j}^{k}}^{+,0}\cup W_{y_{j}^{k}}^{+,1}$, one can show as in
[HS] p 276-8 that if $x\not =y_{l}^{k}$,
\[E_{y_{j}^{k}}^{0}(t)=(\frac{2t}{\pi})^{\frac{n-1}{4}}\mid at
\mid^{1/6}e^{-td(x,y_{j}^{k})}\left(dx_{1}^{(j)}\wedge \ldots
\wedge dx_{k}^{(j)}+O(t^{-1})\right)\]
where $d(x,y_{j}^{k})$ is the Agmon distance between x and $y_{j}^{k}$,
i.e. w.r.t. the metric $\mid df\mid^{2}dg$, but it
is not necessarily true for $x=y_{l}^{k}$.

So we let \[E_{y_{j}^{k}}^{0}(t)(y_{l}^{k})=
(\frac{2t}{\pi})^{\frac{n-1}{4}}\mid at\mid^{1/6}e^{-td(y_{l}^{k},
y_{j}^{k})}c_{lj}(t)\]
where $c_{lj}(t)\in \Lambda^{k}(T_{y_{l}^{k}}(M))$.  By [HS],
\[\mid c_{lj}(t)\mid= O(e^{\epsilon t})\mbox{ for any }\epsilon >0\]

For $x\in U_{y_{l}^{k}}\cap W_{y_{l}^{k}}^{-,0}$,
\[f(x)=f(y_{l}^{k})-(x_{1}^{(l)})^{2}-\ldots-(x_{k}^{(l)})^{2}\]
Therefore,
\[\int_{W_{y_{l}^{k}}^{-,0}}e^{tf}E_{y_{j}^{k}}^{0}(t)=
e^{tf(y_{l}^{k})}\int_{U_{y_{l}^{k}}\cap W_{y_{l}^{k}}^{-,0}}
e^{-t\left((x_{1}^{(l)})^{2}+\ldots +(x_{k}^{(l)})^{2}\right)}
E_{y_{j}^{k}}^{0}(t)+\mbox{ lower order terms}\]
By the stationary phase approximation formula ([D] p23-4), we have
\[ \begin{array}{ll}
    &  \int_{W_{y_{l}^{k}}^{-,0}}e^{tf}E_{y_{j}^{k}}^{0}(t) \\
= & e^{tf(y_{l}^{k})}(\frac{2t}{\pi})^{\frac{n-1-2k}{4}}\mid at
\mid^{\frac{1}{6}}e^{-td(y_{l}^{k},y_{j}^{k})}c_{lj}(t)
\left(\partial_{x_{1}^{(l)}},\ldots,\partial_{x_{k}^{(l)}}\right)+
\mbox{ lower order terms}\\
= & (\frac{2t}{\pi})^{\frac{n-1-2k}{4}}\mid at\mid^{1/6}e^{tf(y_{j}^{k})}
\Xi_{1}(0)\beta_{lj}(t) \end{array} \]
since $d(y_{l}^{k},y_{j}^{k})=f(y_{l}^{k})-f(y_{j}^{k})$ and for some
$\beta_{lj}(t)$.  Hence (ii) is proved.

Note that (cf.[HS] p265, [HS1] p138)
\[\mid \beta_{lj}(t)\mid=O(e^{\epsilon t})\mbox{ for any }\epsilon>0\]
(iii) Since $d(t)E_{y_{j}^{k}}^{0}(t)=\lambda(t)E_{y_{j}^{k}}^{1}(t)$
for some $\lambda(t)\not=0$ when t is sufficiently
large ( this follows from (7) below) and
\[Int_{k+1}\left(e^{tf}d(t)E_{y_{j}^{k}}^{0}(t)\right)=\delta\left
(Int_{k}e^{tf}E_{y_{j}^{k}}^{0}(t)\right) \]
by (ii), we have
\[Int_{k+1}e^{tf}\left(E_{y_{j}^{k}}^{1}\right)=
(\frac{2t}{\pi})^{\frac{n-1-2k}{4}}e^{tf(y_{j}^{k})}\frac{\Xi_{1}(0)
(\mid at\mid )^{1/6}}{\lambda(t)}\left(\delta(e_{y_{j}^{k}}^{0}+
\sum_{l\not = j}\beta_{lj}(t)e_{y_{l}^{k}}^{0})+
O(t^{-1})\right) \]
(iii) is proved by noting that
\[\lambda(t)=\sqrt{e_{1}}\mid at\mid^{1/2}+\mbox{ lower order terms }
\Box \]
Using Proposition 4(i), we can prove

\noindent{\bf Theorem $2''$: (Helffer-Sj\"{o}strand)}  There exist
orthonormal bases $\left\{E_{x_{j}^{k}}(t)\right\}$ of
$\Omega_{small}^{k}(M,t)$, $\left\{E_{y_{j}^{k}}^{0}(t),
E_{y_{j}^{k}}^{1}(t)\right\}$ of $\Omega_{large,k,j}^{\ast}(M,t)$ s.t.
\[<E_{x_{i}^{k+1}}(t),d(t)E_{x_{j}^{k}}(t)>=e^{-t}\left(\sqrt{\frac{t}
{\pi}}\sum_{\gamma}\epsilon_{\gamma}^{new}+
O(t^{-1/2})\right)\]
\[<E_{y_{j}^{k}}^{1}(t),d(t)E_{y_{j}^{k}}^{0}(t)>=
\sqrt{e_{1}}(a_{j}^{(k)}t)^{1/3}+O(t^{1/6})\]
\[<E_{y_{j_{1}}^{k}}^{i_{1}}(t),d(t)E_{y_{j_{2}}^{l}}^{i_{2}}(t)>=
0\mbox{ if $j_{1}\not =j_{2}$ for t sufficiently
large}\]
\[<E_{x_{i}^{k}}(t),d(t)E_{y_{j}^{l}}^{i^{'}}(t)>=
<E_{y_{j}^{l}}^{i^{'}}(t),d(t)E_{x_{i}^{k}}(t)>=0\mbox{ for t
sufficiently large}\]
where $\sum_{\gamma}\epsilon_{\gamma}^{new}=I(x_{i}^{k+1},x_{j}^{k})$
is the incidence number between $x_{i}^{k+1}$ and
$x_{j}^{k}$ defined in \S 1.

\noindent{\bf Proof}: We prove the first and second equalities,
the others are obvious.

Let $d(t)E_{x_{j}^{k}}(t)=\sum_{i}\lambda_{ij}(t)E_{x_{i}^{k+1}}(t)$
for some $\lambda_{ij}(t)$.

Since \[Int_{k+1}e^{tf}\left(d(t)E_{x_{j}^{k}}(t)\right)=
\delta Int_{k}e^{tf}\left(E_{x_{j}^{k}}(t)\right)\]
By Proposition 4(i), we have
\[\begin{array}{ll}
  & \sum_{i}\lambda_{ij}(t)(\frac{2t}{\pi})^{\frac{n-2k-2}{4}}
  e^{t(k+1)}\left(e_{x_{i}^{k+1}}+\sum_{l}I(y_{l}^{k+1},
  x_{i}^{k+1})e_{y_{l}^{k+1}}^{0}+O(t^{-1})\right)\\
= & (\frac{2t}{\pi})^{\frac{n-2k}{4}}\left(\delta e_{x_{j}^{k}}+
\sum_{l}I(y_{l}^{k},x_{j}^{k})\delta e_{y_{l}^{k}}^{0}+
O(t^{-1})\right)
\end{array} \]
By comparing the coefficients of $e_{x_{i}^{k+1}}$,
\[\lambda_{ij}(t)(\frac{2t}{\pi})^{-1/2}e^{t}=i(x_{i}^{k+1},x_{j}^{k})+
\sum_{l}I(y_{l}^{k},x_{j}^{k})i(x_{i}^{k+1},
y_{l}^{k})+O(t^{-1})\]
But by definition of $I(x,y)$,
\[I(x_{i}^{k+1},x_{j}^{k})=i(x_{i}^{k+1},x_{j}^{k})+\sum_{l}I(y_{l}^{k},
x_{j}^{k})i(x_{i}^{k+1},y_{l}^{k})=
\sum_{\gamma}\epsilon_{\gamma}^{new}\]
Hence the first equality is proved. For the second equality, note that
$E_{y_{j}^{k}}^{i}(t)$ are normalized eigenforms
in $\Omega_{large,k,j}^{\ast}(M,t)$ and
\begin{equation}
\|d(t)E_{y_{j}^{k}}^{0}(t)\|^{2}=<e_{y_{j}^{k}}^{0}(t),\Delta^{k}(t)
E_{y_{j}^{k}}^{0}(t)>=e_{1}\mid a_{j}^{(k)}t
\mid^{2/3}+\mbox{ lower order terms}
\end{equation}
So,

\[<e_{y_{j}^{k}}^{1}(t),d(t)E_{y_{j}^{k}}^{0}(t)>=\pm\sqrt{e_{1}}
(\mid a_{j}^{(k)}\mid t)^{1/3}+\mbox{ lower order
terms} \mbox{  }\Box \]
In view of Proposition 4, define $f^{k}(t):\Omega_{0}^{k}(M,t)
\rightarrow C^{k}(M,f)$ s.t.
\[\begin{array}{lll}
   f^{k}(t)\left(E_{x_{j}^{k}}(t)\right) & = &
   (\frac{\pi}{2t})^{\frac{n-2k}{4}}e^{-tk}Int_{k}e^{tf}
   \left(E_{x_{j}^{k}}(t)\right)  \\
   f^{k}(t)\left(E_{y_{j}^{k}}^{0}(t)\right) & = & Int_{k}e^{tf}\left
   (E_{y_{j}^{k}}^{0}(t)\right)  \\
   f^{k+1}(t)\left(E_{y_{j}^{k}}^{1}(t)\right) & = & Int_{k+1}e^{tf}
   \left(E_{y_{j}^{k}}^{1}(t)\right)
  \end{array}   \]
Let \[\begin{array}{lll}
       \left(\Omega_{0}^{\ast}(M,t),\tilde{d}(t)\right) & = &
	   \left(\Omega_{small}^{\ast}(M,t),
	   e^{t}(\pi/2t)^{1/2}d(t)\right) \\
                                                        &   & \perp\left
														(\perp_{k,j}\left
				(\Omega_{large,k,j}^{\ast}(M,t),d(t)\right) \right)
      \end{array}   \]
Also define
\[\begin{array}{l}
\hat{e}_{x_{j}^{k}}=e_{x_{j}^{k}}+\sum_{l}I(y_{l}^{k},x_{j}^{k})
e_{y_{l}^{k}}^{0}\\
\hat{e}_{y_{j}^{k}}^{0}=e_{y_{j}^{k}}^{0}\\
\hat{e}_{y_{j}^{k}}^{1}=\delta(e_{y_{j}^{k}}^{0})
\end{array} \]
Then we have

\begin{proposition}  $f^{\ast}(t):(\Omega_{0}^{\ast}(M,t),\tilde{d}(t))
\rightarrow (C^{\ast}(M,f),\delta)$
is a morphism of co-chain complexes s.t.
\[f^{k}(t)\left(E_{x_{j}^{k}}(t)\right)=\hat{e}_{x_{j}^{k}}+O(t^{-1})  \]
with $f^{k}(t)\left(E_{y_{j}^{k}}^{0}(t)\right),f^{k+1}
\left(E_{y_{j}^{k}}^{1}(t)\right)$ given by (ii) and (iii)
in Proposition 4.
\end{proposition}

Let the matrix associated to the linear map $f^{k}(t)$ w.r.t. the bases
\[\left\{E_{x_{j}^{k}}(t),E_{y_{j}^{k}}^{0}(t),E_{y_{j}^{k-1}}^{1}(t)
\right\}\mbox{ and } \left\{\hat{e}_{x_{j}^{k}},
e_{y_{j}^{k}}^{0},e_{y_{j}^{k-1}}^{1}\right\}\]
be
\[F^{k}(t)=\left(\begin{array}{ccc}
I & O(e^{kt}) & N_{1}^{k}(t) \\
O(t^{-1}) & M^{k}(t) & N_{2}^{k}(t) \\
O(t^{-1}) & O(e^{kt}) & N_{3}^{k}(t)
\end{array}  \right)   \]
where $O(t^{-1})$ in a certain entry of the matrix means that the
corresponding entry is of the order $O(t^{-1})$.
Here $M^{k}(t)$ is $\left(\frac{2t}{\pi}\right)^{\frac{n-1-2k}{4}}
\Xi_{1}(0)$ times the following matrix
\[\left(\begin{array}{cccc}

\mid a_{1}^{(k)}t\mid^{\frac{1}{6}}e^{tf(y_{1}^{k})} & \mid a_{2}^{(k)}
t\mid^{\frac{1}{6}}e^{t(f(y_{1}^{k})-
\epsilon_{0})}\beta_{12}(t) & \cdots & \mid a_{m_{k}^{'}}^
{(k)}t\mid^{\frac{1}{6}}e^{t(f(y_{1}^{k})-\epsilon_{0})}
\beta_{1m_{k}^{'}}(t) \\
\mid a_{1}^{(k)}t\mid^{\frac{1}{6}}e^{tf(y_{1}^{k})}\beta_{21}(t) &
\mid a_{2}^{(k)}t\mid^{\frac{1}{6}}e^{tf(y_{2}^{k})} &
\cdots & \mid a_{m_{k}^{'}}^{(k)}t\mid^{\frac{1}{6}}e^{t(f(y_{2}^{k})-
\epsilon_{0})}\beta_{2m_{k}^{'}}(t) \\
\vdots & \vdots & \ddots & \vdots \\
\mid a_{1}^{(k)}t\mid^{\frac{1}{6}}e^{tf(y_{1}^{k})}
\beta_{m_{k}^{'}1}(t) & \mid a_{2}^{(k)}t\mid^{\frac{1}{6}}
e^{tf(y_{2}^{k})}\beta_{m_{k}^{'}2}(t) & \cdots & \mid
a_{m_{k}^{'}}^{(k)}t\mid^{\frac{1}{6}}e^{tf(y_{m_{k}^{'}}^{k})}
\end{array} \right) \]
Note that we have used the fact that the birth-death points are
indexed such that
\[f(y_{1}^{k})\leq f(y_{2}^{k})\leq \ldots \leq f(y_{m_{k}^{'}}^{k})\]
so that the above matrix is approximately lower triangular. Hence,
$M^{k}(t)$ is invertible for sufficiently large t.

Also, let
\[A^{k}(t)=diag\left((\sqrt{e_{1}}\mid a_{1}^{(k)}t\mid^{1/3})^{-1},
\ldots,(\sqrt{e_{1}}\mid a_{m_{k}^{'}}^{(k)}t
\mid^{1/3})^{-1}\right)\]
Define for $1\leq j\leq m_{k}^{'}$,
\[\hat{E}_{y_{j}^{k}}^{0}(t)=\left(M^{k}(t)\right)^{-1}\left
(E_{y_{j}^{k}}^{0}(t)\right)\]
\[\hat{E}_{y_{j}^{k}}^{1}(t)=\left(M^{k}(t)A^{k}(t)\right)^{-1}
E_{y_{j}^{k}}^{1}(t)\]
Observe that $\left\{\hat{E}_{y_{j}^{k}}^{i}(t)\right\}$ is
approximately orthogonal whose elements are still localized
at the corresponding birth-death points.

Let \[B^{k}(t)=\left(\begin{array}{ccc}
I & 0 & 0 \\
0 & \left(M^{k}(t)\right)^{-1} & 0 \\
0 & 0 & \left(M^{k-1}(t)A^{k-1}(t)\right)^{-1}
\end{array} \right) \]
Then the matrix associated to $f^{k}(t)$ w.r.t. the new bases
\[\left\{E_{x_{j}^{k}}(t),\hat{E}_{y_{j}^{k}}^{0}(t),
\hat{E}_{y_{j}^{k}}^{1}(t)\right\}\mbox{ and }
\left\{\hat{e}_{x_{j}^{k}},e_{y_{j}^{k}}^{0},e_{y_{j}^{k-1}}^{1}
\right\} \]
is \[F^{k}(t)B^{k}(t)=\left(\begin{array}{ccc}
I & O(e^{kt}) & N_{1}^{k}(t) \\
O(t^{-1}) & M^{k}(t) & N_{2}^{k}(t) \\
O(t^{-1}) & O(e^{kt}) & N_{3}^{k}(t)
\end{array} \right) \left(\begin{array}{ccc}
I & 0 & 0 \\
0 & \left(M^{k}(t)\right)^{-1} & 0 \\
0 & 0 & \left(M^{k-1}(t)A^{k-1}(t)\right)^{-1}
\end{array} \right) \]
Suppose that \[\delta^{(k-1)}=\left(\begin{array}{ccc}
\delta_{11}^{(k-1)} & \delta_{12}^{(k-1)} & \delta_{13}^{(k-1)} \\
\delta_{21}^{(k-1)} & \delta_{22}^{(k-1)} & \delta_{23}^{(k-1)} \\
\delta_{31}^{(k-1)} & \delta_{32}^{(k-1)} & \delta_{33}^{(k-1)}
\end{array} \right) \]
w.r.t. the bases $\left\{\hat{e}_{x_{j}^{k}},e_{y_{j}^{k}}^{0},
e_{y_{j}^{k}}^{1}\right\}$.

Then, \[F^{k}(t)B^{k}(t)=\left(\begin{array}{ccc}
I & O(t^{-1}) & \delta_{13}^{(k-1)}+O(t^{-1}) \\
O(t^{-1}) & I & \delta_{23}^{(k-1)}+O(t^{-1}) \\
O(t^{-1}) & O(t^{-1}) & \delta_{33}^{(k-1)}+O(t^{-1})
\end{array} \right) \]
To see this, it suffices to show
\begin{equation}
Int_{k}e^{tf}\left(\hat{E}_{y_{j}^{k-1}}^{1}(t)\right)=\delta
(e_{y_{j}^{k-1}}^{0})+O(t^{-1})
\end{equation}
But by Proposition 4(iii),
\[Int_{k}e^{tf}\left(E_{y_{j}^{k-1}}^{1}\right)=\delta\left(\sum_{l}
\left(M^{k-1}(t)A^{k-1}(t)\right)_{lj}
e_{y_{l}^{k-1}}^{0}+O(t^{-1})\right)\]
Using the definition of $\hat{E}_{y_{j}^{k-1}}^{1}(t)$, (8) follows.

Hence, with the definition of $\hat{e}_{y_{j}^{k}}^{i}$ on p18, we
finally have

\noindent{\bf Theorem 2}: $f^{\ast}(t):(\Omega_{0}^{\ast}(M,t),
\tilde{d}(t))\longrightarrow(C^{\ast}(M,f),\delta)$
is a morphism of cochain complexes s.t.
\[f^{\ast}(t)=I+O(t^{-1})\]
w.r.t. the bases $\left\{E_{x_{j}^{k}}(t),\hat{E}_{y_{j}^{k}}^{i}(t)
\right\}$ and $\left\{\hat{e}_{x_{j}^{k}},
\hat{e}_{y_{j}^{k}}^{i}\right\}$.

Inside $(C^{\ast}(M,f),\delta)$, there is a subcomplex
$(C_{nd}^{\ast}(M,f),\delta)$ such that
\[dimC_{nd}^{k}(M,f)=m_{k}\]
where $m_{k}$ is the number of non-degenerate critical points of
index k. This subcomplex can be obtained by
application of the following Lemma.
\begin{lemma} Suppose $(C^{\ast},\delta)$ is a cochain complex such that
\[C^{\ast}=\left\{ \begin{array}{ll}
\tilde{C}^{\ast} & \mbox{ if $\ast\not =k,k+1$}\\
\tilde{C}^{\ast}\oplus R  & \mbox{ if $\ast=k$ or $k+1$}
\end{array} \right. \]
Let $dim(\tilde{C}^{q})=n_{q}$, for $0\leq q\leq n$,
\[\left\{e_{x_{1}^{q}},\ldots,e_{x_{n_{q}}^{q}}\right\} \mbox{ be a
basis of $\tilde{C}^{q}$}\]
so that \[\left\{e_{x_{1}^{k}},\ldots,e_{x_{n_{k}}^{k}},e_{y}^{0}
\right\} \mbox{ is a basis of
$\tilde{C}^{k}\oplus R$}\]
and \[\left\{e_{x_{1}^{k+1}},\ldots,e_{x_{n_{k+1}}^{k+1}},e_{y}^{1}
\right\} \mbox{ is a basis of
$\tilde{C}^{k+1}\oplus R$}\]
and w.r.t. the above bases,

\[\delta^{(k)}=\left( \begin{array}{cccc}
i(x_{1}^{k+1},x_{1}^{k}) & \cdots & i(x_{1}^{k+1},x_{n_{k}}^{k}) &
i(x_{1}^{k+1},y) \\
\vdots & \ddots & \vdots & \vdots \\
i(x_{n_{k+1}}^{k+1},x_{1}^{k}) & \cdots & i(x_{n_{k+1}}^{k+1},
x_{n_{k}}^{k}) & i(x_{n_{k+1}}^{k+1},y) \\
i(y,x_{1}^{k}) & \cdots & i(y,x_{n_{k}}^{k}) & 1 \end{array} \right) \]
Then with the following change of bases in $C^{k}$ and $C^{k+1}$,
\[\left\{e_{x_{l}^{k}},e_{y}^{0}\right\}_{1\leq l\leq n_{k}}\longmapsto
\left\{e_{x_{l}^{k}}-i(y,x_{l}^{k})
e_{y}^{0},e_{y}^{0}\right\}_{1\leq l\leq n_{k}}\]
\[\left\{e_{x_{l}^{k+1}},e_{y}^{1}\right\}_{1\leq l\leq n_{k+1}}
\longmapsto \left\{e_{x_{l}^{k+1}},\delta
e_{y}^{0}=e_{y}^{1}+\sum_{j=1}^{n_{k+1}}i(x_{j}^{k+1},y)e_{x_{j}^{k+1}}
\right\}_{1\leq l\leq n_{k+1}} \]
we have

(i)\[\delta^{(k)}=\left(\begin{array}{cccc}
i^{'}(x_{1}^{k+1},x_{1}^{k}) & \cdots & i^{'}(x_{1}^{k+1},
x_{n_{k}}^{k}) & 0 \\
\vdots & \ddots & \vdots & \vdots \\
i^{'}(x_{n_{k+1}}^{k+1},x_{1}^{k}) & \cdots & i^{'}(x_{n_{k+1}}^{k+1},
x_{n_{k}}^{k}) & 0 \\
0 & \cdots & 0 & 1 \end{array} \right) \]
where $i^{'}(x_{i}^{k+1},x_{j}^{k})=i(x_{i}^{k+1},x_{j}^{k})-
i(x_{i}^{k+1},y)i(y,x_{j}^{k})$.

(ii) \[\delta^{(k-1)}=\left( \begin{array}{ccc}
i(x_{1}^{k},x_{1}^{k-1}) & \cdots & i(x_{1}^{k},x_{n_{k-1}}^{k-1}) \\
\vdots & \ddots & \vdots \\
i(x_{n_{k}}^{k},x_{1}^{k-1}) & \cdots & i(x_{n_{k}}^{k},
x_{n_{k-1}}^{k-1}) \\
0 & \cdots & 0 \end{array} \right) \]
\[\delta^{(k+1)}=\left(\begin{array}{cccc}
i(x_{1}^{k+2},x_{1}^{k+1}) & \cdots & i(x_{1}^{k+2},
x_{n_{k+1}}^{k+1}) & 0 \\
\vdots & \ddots & \vdots & \vdots \\
i(x_{n_{k+2}}^{k+2},x_{1}^{k+1}) & \cdots & i(x_{n_{k+2}}^{k+2},
x_{n_{k+1}}^{k+1}) & 0 \end{array} \right) \]
\end{lemma}
{\bf Corollary}: Let
\[(C^{'})^{\ast}=\left\{ \begin{array}{ll}
C^{\ast} & \mbox{ if $\ast\not =k$} \\
Span\{e_{x_{l}^{k}}-i(y,x_{l}^{k})e_{y}^{0}\}_{1\leq l\leq n_{k}} &
\mbox{ if $\ast=k$} \end{array} \right.\]
\[(C^{''})^{\ast}=\left\{ \begin{array}{ll}
0 & \mbox{ if $\ast\not =k,k+1$} \\
e_{y}^{0} & \mbox{ if $\ast=k$} \\
\delta e_{y}^{0} & \mbox{ if $\ast=k+1$} \end{array} \right. \]
Then \[(C^{\ast},\delta)=\left((C^{'})^{\ast},\delta\right)\oplus
\left((C^{''})^{\ast},\delta\right)\]
In particular, $\left((C^{'})^{\ast},\delta\right)$ is a subcomplex of
$(C^{\ast},\delta)$. Both of them calculate
the same cohomology.

\noindent{\bf Remark}: For the application of Lemma 3, it is clear
that $\tilde{C}^{\ast}$ need not be generated only
by cells corresponding to non-degenerate critical points, but it can
also be generated by cells corresponding to
birth-death points.  Let
\[k<f(y_{1}^{k})\leq \ldots \leq f(y_{m_{k}^{'}}^{k})<k+1\]
We $'$eliminate$'$ $y_{1}^{k}$ first by applying Lemma 3 and obtain
a subcomplex $\left((C^{(1)})^{\ast},\delta\right)$.
Note that \[\left\{ \begin{array}{l}
e_{y_{2}^{k}}^{0}=e_{y_{2}^{k}}^{0}-i(y_{1}^{k},y_{2}^{k})
e_{y_{1}^{k}}^{0}\in ((C^{(1)})^{k}\\
e_{y_{2}^{k}}^{1}\in (C^{(1)})^{k+1} \\
\delta(e_{y_{2}^{k}}^{0})=e_{y_{2}^{k}}^{1}+\ldots \end{array}
\right. \]
Therefore, the assumptions in Lemma 3 are satisfied and we can apply
Lemma 3 to $'$eliminate$'$ $y_{2}^{k}$ and obtain
$\left((C^{(2)})^{\ast},\delta\right)$.  Hence by applying Lemma 3
repeatedly, $\left(C_{nd}^{\ast}(M,f),\delta\right)$
is obtained.

\noindent{\bf Proof of Lemma 3}:  The Lemma can be proved by
direct calculation.

\noindent{\bf Theorem $2''$}:  $f^{\ast}(t)\mid_{\Omega_{small}
^{\ast}(M,t)}:\left(\Omega_{small}^{\ast}(M,t),\tilde{d}(t)
\right)\longrightarrow(C^{\ast}(M,f),\delta)$ is an injective
homomorphism of cochain complexes whose image complex
converges to $(C_{nd}^{\ast}(M,f),\delta)$ in $(C^{\ast}(M,f),
\delta)$ as $t\rightarrow\infty$, more precisely,
\[f^{k}(f)\left(E_{x_{j}^{k}}(t)\right)=\hat{e}_{x_{j}^{k}}+O(t^{-1})
\mbox{ in }C^{\ast}(M,f)\]
{\bf Remark}: One can show by induction that
\[\hat{e}_{x_{j}^{k}}=e_{x_{j}^{k}}+\sum_{l}I(y_{l}^{k},x_{j}^{k})
e_{y_{l}^{k}}^{0}\in C_{nd}^{k}(M,f)\]

\noindent{\Large\bf Appendix}

\noindent{\bf Sketch of Proof (of Theorem 1)}:

Let $C_{bd}$ be the set of birth-death critical points of f,

 $C_{nd}$ be the set of non-degenerate critical points of f,

 $0=e_{1}^{k}= \ldots =e_{m_{k}}^{k}$ be the smallest eigenvalues of
 $\bigoplus _{j\in C_{nd}}\overline{\Delta_{j}^{k}}(1)$

 $0\leq e_{m_{k}+1}^{k} \leq e_{m_{k}+2}^{k} \leq \ldots$ be the
 eigenvalues of $\bigoplus_{j\in C_{bd}}
 \overline{\Delta_{j}^{k}}(1)$

Let $\left\{\Psi_{l}^{k}(1)\right\}_{l}^{\infty}$ be the eigenvectors
corresponding to $e_{l}^{k}$ of
$\bigoplus_{j\in Crit(f)}\overline{\Delta_{j}^{k}}(1)$

More generally, let $\left\{\Psi_{l}^{k}(t)\right\}_{l}^{\infty}$ be
the eigenvectors corresponding to $e_{l}^{k}t^{2/3}$
of $\bigoplus_{j\in Crit(f)}\overline{\Delta_{j}^{k}}(t)$.

Note that \[ \Psi_{l}^{k}(t) = \left\{ \begin{array}{ll}
                                       U(t^{1/2})\Psi_{l}^{k}(1) &
									   \mbox{if $1\leq l \leq m_{k}$} \\
                                       U(t^{1/3})\Psi_{l}^{k}(1) &
									   \mbox{if $m_{k}+1\leq l$}
                                       \end{array}
                               \right.   \]

Then Theorem 1 is essentially equivalent to \[\lim_{t\rightarrow \infty}
\frac{E_{l}(t)}{t^{2/3}}=e_{l}^{k} \]
The proof is divided into 2 steps.

(i) $\overline{\lim}_{t \rightarrow \infty}\frac{E_{l}(t)}{t^{2/3}}
\leq e_{l}^{k}$

This follows by similar arguments as in [S].  Let $\left\{J_{j}
\right\}_{j\in Crit(f)\cup\{0\}}$ be a partition of
unity on M (cf.[S] p 27).  Let $\Psi_{l}^{k}(t)$ be an eigenvector of
$\overline{\Delta_{j(l)}^{k}}(t)$, define
$\varphi_{l}^{k}(t)=J_{j(l)}\Psi_{l}^{k}(t)$. Then
$\left\{\varphi_{l}^{k}(t)\right\}$ form a set of approximate
eigenvectors in $L^{2}(\Lambda^{k}(M))$ with
\[<\varphi_{m}(t),\Delta^{k}(t)\varphi_{n}(t)>=e_{n}^{k}t^{2/3}
\delta_{nm}+o(t^{2/3}) \]
by using the minimax principle, (i) follows.

(ii) $\underline{\lim}_{t \rightarrow \infty}\frac{E_{l}(t)}{t^{2/3}}
\geq e_{l}^{k}$

To prove (ii), one has to modify slightly the arguments in [S].  It
suffices to construct, for $e\in(e_{l}^{k},
e_{l+1}^{k})$, a symmetric operator $R(t)$ of rank $l$ s.t.
\[\Delta^{k}(t) \geq t^{2/3}e+R(t)+o(t^{2/3}) \ldots(\ast) \]

To construct $R(t)$ define $\Delta_{j}^{k}(t):C^{\infty}(\Lambda^{k}(M))
\rightarrow C^{\infty}(\Lambda^{k}(M))$
\[\Delta_{j}^{k}(t)=\Delta^{k}+t^{2}f_{j}+tA \]
\[f_{j}(x) = \left\{ \begin{array}{ll}
                      \mid df(x) \mid^{2} & \mbox{if $x\in U_{j}$}  \\
                      >0                  & \mbox{if $x \not \in U_{j}$}
                     \end{array}
             \right.   \]
Let $0\leq E_{1}^{(j)}(t)\leq E_{2}^{(j)}(t)\leq \ldots \leq
E_{l}^{(j)}(t)\leq \ldots $ be the eigenvalues of
$\Delta_{j}^{k}(t)$,

$\Psi_{1}^{(j)}(t),\psi_{2}^{(j)}(t),\ldots,\Psi_{l}^{(j)}(t),\ldots $
be the corresponding eigenvectors of
$\Delta_{j}^{k}(t)$.

\noindent For $j\in C_{bd}$, let $0\leq e_{1}^{(j)}\leq e_{2}^{(j)}\leq
\ldots \leq e_{j}^{(j)}\leq \ldots $ be the
eigenvalues of $\overline{\Delta_{j}^{k}}(1)$ ,then one can show that

\[\lim_{t\rightarrow \infty}
\frac{E_{l}^{(j)}}{t^{2/3}}=e_{l}^{(j)} \]
Define $n_{j}$ s.t. $e_{n_{j}}^{(j)}<e<e_{n_{j}+1}^{(j)}$
\[P_{j}(t)=\left\{ \begin{array}{ll}
                   \mbox{orthogonal projection onto $span\left\{
				   \Psi_{l}^{(j)}\right\}_{1\leq l\leq n_{j}}$} &
				   \mbox{if $j\in C_{bd}$}  \\
                   \mbox{orthogonal projection onto smallest eigenvector
				   of $\Delta_{j}^{k}(t)$} & \mbox{if $j\in C_{nd}$}
                   \end{array}
           \right.  \]
$R_{j}(t)=\left(\Delta_{j}^{k}(t)-t^{2/3}e\right)P_{j}(t)$

Define \[R(t)=\sum_{j\in Crit(f)}J_{j}R_{j}(t)J_{j}\]

which is a symmetric operator of rank $l$.

To verify $(\ast)$, observe that by IMS localization formula (cf.[S] p28)

\[\Delta^{k}(t)\geq t^{2/3}eJ_{0}^{2}+\sum_{j\in Crit(f)}J_{j}
\Delta_{j}^{k}(t)J_{j}+O(1)\]
Then $(\ast)$ follows from the definition of $R(t)$.

One finishes the proof by showing that
\[\lim_{t\rightarrow \infty}E_{1}(t)= \ldots =\lim_{t\rightarrow
\infty}E_{m_{k}}(t)=0
\mbox{ $$ $$ $$ $$ $$  $\Box$}  \]

\end{document}